\magnification=\magstep1
\noindent

\input amstex
\UseAMSsymbols
\input pictex 
\vsize=23truecm

\NoBlackBoxes
\parindent=18pt
  
   \font\rmk=cmr8      
   
\font\gross=cmbx10 scaled\magstep1

\def\mod{\operatorname{mod}}

\def\Hom{\operatorname{Hom}}
\def\End{\operatorname{End}}
\def\Ext{\operatorname{Ext}}

\def\rad{\operatorname{rad}}
\def\add{\operatorname{add}}

\def\soc{\operatorname{soc}}

\def\D{\operatorname{D}}

\def\ss{\ssize}

\def\arr#1#2{\arrow <1.5mm> [0.25,0.75] from #1 to #2}

\def\s{\hfill \square}


\vglue1truecm	

\centerline{\gross The brick chain complexity of an artin algebra.}
	\medskip
\centerline{Claus Michael Ringel}
	\bigskip\medskip
{\narrower {\bf Abstract.} We consider the category of finitely generated modules over an
artin algebra $A$. It is known that any module $M$ has a 
brick chain filtration. We say that $M$ has brick chain complexity 
at most $t$ provided $M$ has a brick chain filtration of length at most $t$. The brick chain
complexity of $A$ is by definition the supremum of the brick chain complexity of
the indecomposable $A$-modules.
The aim of this note is to calculate the brick chain complexity for some algebras.
We will exhibit algebras with arbitrarily large brick chain complexity.
\par} 
	\bigskip\medskip
{\bf 1. Brick chain filtrations.}
	\medskip
{\bf 1.1.} We deal with an artin algebra $A$; the modules to be considered are the
left $A$-modules of finite length. 
Given a set $\Cal X$ of modules, let $\Cal E(\Cal X)$ be the class of modules
which have a filtration with all factors in $\Cal X$. 	
We recall that a {\it brick} is a module whose endomorphism ring is a division ring.
If $B$ is a brick, the modules in $\Cal E(B)$ will be said to be {\it homogeneous}
of type $B$.
A finite sequence $(B_1,\dots, B_m)$ of bricks is called 
a {\it brick chain,} provided $\Hom(B_i,B_j) = 0$ for $i<j$. 
A filtration $0 = M_0 \subset M_1 \subset \cdots \subset M_m = M$ is called
a {\it brick chain filtration} of $M,$ provided there is a brick chain
$(B_1,\dots, B_m)$ (its {\it type}) such that $M_i/M_{i-1}$ is homogeneous of type $B_i$,
for all $1 \le i \le m$; the number $m$ is called the {\it length} of the filtration. 
We have shown in [R1] that {\it any module has at least one (and usually several) brick chain
filtrations.} 
	\medskip
{\bf 1.2. Numerical data provided by brick chain filtrations.} 
The consideration of brick chain filtrations $(M_i)_i$ 
of a module $M$ should help to separate different features of $M$; 
the various data provided by any brick chain filtration sheds light
on the structure of the module. 
Also, if there are given several brick chain filtrations of the same module, one may ask
in which way the respective data are related to each other. 
	\medskip
At first sight, there are at least four different directions for investigations.
	\smallskip
 \item{$\bullet$} First of all, there is the number of bricks in the category $\mod A$: 
   is it finite or infinite? Brick finiteness is known as an important finiteness condition.  
   Given any mdoule $M$, one should look at the set of bricks which occur in the brick type
   of a brick chain filtration of $M$ (see 1.4); it is an open question whether this set
   is always finite (see [R1], 3.2). 
 \item{$\bullet$} Second, given a module $M$, one may look at 
   the minimal length of a brick chain filtration of $M$, we call it the {\it brick chain
   complexity} of $M$. 
   In the present note we look at the brick chain complexity 
   of indecomposable modules $M$.
 \item{$\bullet$}	
   Third, given a brick $B$, the Loewy length of the category $\Cal E(B)$
   seems to play an important role, we call it the {\it Loewy number} $\lambda(B)$
   of $B$. Note that the Loewy number of a brick is a positive integer or $\infty.$
   We will deal with the Loewy number of some bricks in [R2]. 
 \item{$\bullet$}	
   Forth, given a brick $B$, we denote by $e(B)$ the {\it e-dimension} of $B$,
   by definition this is the dimension of $\Ext^1(B,B)$ as an
   $\End(B)$-vector space. The "e" refers to "extension", or also to "embedding"
   (in commutative algebra, given a local ring $R$, the length of
   $\rad R/\rad^2R$ as an $R$-module is called the embedding dimension).
   The e-dimension of a brick is always a non-negative integer.   
	\smallskip
Altogether, we note that the various data mentioned here provide many different ways for
specifying properties of modules and algebras. 
Of course, these numerical values are only a bare shadow
of the information which is available. For example, instead of looking at the Loewy number
and the e-dimension of $B$ one should take into account the complete structure
of $\Cal E(B)$.
	\medskip
{\bf 1.3. Building blocks of modules.}
As the wording "brick" indicates, we envision bricks as building blocks for modules.
Given a module $M$, there are the various bricks which occur in the brick type 
of a brick chain filtration of $M$; these are the bricks which should be considered as 
the building blocks for $M$. (Why do we want to take into account all brick chain
filtrations, not just the torsional ones? Of course, the torsional ones are
well behaved and easy to overlook; however, as we know, a brick chain filtration
induced by a torsional one usually is not torsional: this seems to indicate that the 
class of torsional brick chain filtrations is too restrictive.) 

	\bigskip\medskip
{\bf 2. Brick chain complexity. Definition and some examples.}
	\medskip
{\bf 2.1.} As we have mentioned, a module $M$ is said to have {\it 
brick chain complexity} at most $t$ provided there is a brick chain filtration 
with at most $t$ factors. The {\it brick chain complexity} of an algebra $A$ is the 
supremum of the brick chain complexity of the indecomposable $A$-modules
(it is a natural number or $\infty$). 
	\medskip
{\bf 2.2. Algebras with brick chain complexity 1.}
Of course, bricks are modules with brick chain complexity 1. Thus, any non-zero algebra
which is representation-directed has complexity 1. Since any indecomposable
Kronecker module is homogeneous, the Kronecker algebra has also complexity 1.
Next, if $A$ is a local algebra, then again all modules are
homogeneous, thus local algebras also have complexity 1.
	\medskip
{\bf 2.3. Algebras with brick chain complexity 2.}
Nakayama algebras, tame concealed algebras and tubular algebras have
brick chain complexity at most 2. For example, 
if $A$ is tame concealed, the only indecomposable modules which are not
homogeneous are the indecomposable modules $M$ which belong to a tube say of rank $r$,
with regular length not divisible by $r,$ and these modules have complexity 2.
	\medskip
{\bf 2.4. An example of modules with brick chain complexity 3.}
Note that modules which belong to a standard tube 
may have brick chain complexity 3. As an example, take the non-stable tube for the algebra
$$
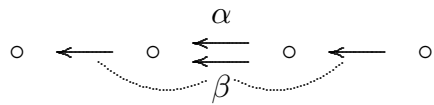
{\beginpicture
    \setcoordinatesystem units <1.5cm,1cm>
\put{} at 0 -.5
\multiput{$\circ$} at 0 0  1 0  2 0  3 0 /
\arr{0.7 0}{0.3 0}
\arr{1.7 0.1}{1.3 0.1}
\arr{1.7 -.1}{1.3 -.1}
\arr{2.7 0}{2.3 0}
\setdots <.4mm>
\setquadratic
\plot 0.6 -.1  1 -.4  1.4 -.3 /
\plot 1.6 -.3  2 -.4  2.4 -.1 /
\put{$\alpha$} at 1.5 0.4
\put{$\beta$} at  1.5 -.4
\endpicture}
$$
{\it An indecomposable module $M = M(t)$ with dimension vector $(1,t,t,1)$ has
brick chain complexity $\min(t,3)$.} Let us exhbit the module $M(4)$ explicitly:
the action of $\beta$ is given by dashed lines, the action of all other arrows by 
solid lines. The coefficient quiver of $M(4)$ is the following (with arrows going
from right to left):
$$
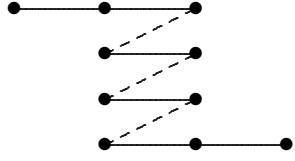
{\beginpicture
    \setcoordinatesystem units <1cm,.5cm>
\multiput{} at 0 0  3 3 /
\multiput{$\bullet$} at 0 3  1 3  2 3  1 0  1 1  1 2  2 0  2 1  2 2  3 0 /
\plot 0 3  2 3 /
\plot 1 2  2 2 /
\plot 1 1  2 1 /
\plot 1 0  3 0 /
\setdashes <1mm> 
\plot 1 2  2 3 /
\plot 1 1  2 2 /
\plot 1 0  2 1 /
\endpicture}
$$
	\medskip

{\bf 2.5. Algebras with finite brick chain complexity.}
If the quiver of the algebra $A$ is directed, then the brick chain complexity
is bounded by the Loewy length of $A$. In this way, we see that the brick chain complexity
of the Kronecker algebra is at most $2$, but this bound is not optimal: 
as we have mentioned already, the brick chain complexity of the Kronecker algebra is
equal to $1$.
	\medskip
{\bf 2.6. Target of this note.} 
The target of the present note is to present algebras with arbitrary 
brick chain complexity. 
This will be done in sections 4 and 6: In section 4 we present algebras with finite
brick chain complexity, in section 6 algebras with infinite brick chain complexity.
In section 3 we study properties of the bricks $B$
which occur as the first entry in the brick type of a brick chain filtration of a module $M$, 
we call them
``foundation bricks'' of $M$. Section 5 provides
some observations abour Kronecker modules which will be used in section 6.
 
Let us stress that it should not be considered as surprising 
that algebras can have infinite brick chain complexity:
all sufficiently complicated algebras will have this property. However, given 
a brick chain filtration of a module, usually it is not easy to decide whether
the length of this filtration is smallest possible, thus the problem is just
to verify that a given algebra has infinite brick chain complexity. 
The algebra which we present has the advantage 
that there is a class of modules where it is rather easy to observe that 
all brick chain filtrations have the same length. 
	\bigskip\medskip
{\bf 3. Foundation bricks.}
	\medskip
A brick $B$ is said to be a {\it foundation brick} for the module $M$ provided 
there is a brick chain filtration of $M$ of type $(B_1,B_2,\dots,B_m)$ with $B_1 = B.$
	\medskip
{\bf 3.1. Lemma.}
{\it Let $B$ be a foundation brick for $M$. 
Then the sum $t_BM$ of the images $B \to M$ is a direct sum of copies of $B.$}
	\medskip
We recall that the sum $t_B(M)$ of the images of the maps $B \to M$ is 
called the {\it trace} of $B$ in $M$.
	\medskip
Proof: Since $B$ is a foundation brick for $M$, 
there is a brick chain filtration $(M_i)_i$ of $M$ of type $(B_1,\dots,B_m)$
with $B = B_1.$ Let $f\:B \to M$ be a homomorphism. Since $\Hom(B,M/M_1) = 0$, we see that
$f$ maps into $M_1,$ thus $f$ is a map inside of $\Cal E(B).$  
This shows that the trace of $B$ in $M$ is just the
relative socle of $M_1$ as an object in $\Cal E(B),$ 
thus it is a direct sum of copies of $B$.
$\s$
	\medskip
{\bf Remark.} Let us stress that the converse of 3.1 is not true: Let $M$ be a serial
module with composition factors $1,1,2$ (going upwarrds), and $B$ its submodule of
length $2$. Then $B$ is not a foundation brick of $M$, however the trace of $M$ in $M$ 
is (of course) a direct sum of copies of $B$.
	\medskip 
{\bf 3.2. Corollary.} {\it Let $B$ be a foundation brick of the module $M$.
If a simple factor module of $B$ is isomorphic to a submodule of $M$, then $B$ is simple.}
	\medskip
Proof. Let $U$ be the trace of $B$ in $M$. According to the Lemma, $U = \bigoplus_i U_i$
for some submodules $U_i$ of $M$ isomorphic to $B$. Let $f\:B \to M$ be a map with simple
image $S$. Then $S$ is a submodule of the trace $U$ of $B$ in $M$. It follows that $S$ is
a submodule of some $U_i$, thus a submodule of $B$. But $S$ is also a factor module of $B$.
Since $B$ is a brick, it follows that $B = S.$
$\s$
	\medskip
{\bf Remark.} By definition, the {\it socle bricks} of a module $M$ are the indecomposable
direct summands of the iterated endo-socle of $M$. 
{\it Let $M$ be a module. Any socle brick of $M$ is a foundation brick of
$M$, but the converse does not hold.} 

Proof. Let $B$ be a socle brick of $M$. Then $\D B$
is a top brick of $N = \D M$, thus there is a brick chain filtration $(N_i)_i$ of $N$
of some type $(C_1,\dots,C_m)$ with $C_m = \D B.$ The duality $\D$ yields a
brick chain filtration $(M_i)_i$ of $M$ with 
$N_i = \D N/\D N_{m-i}$, its type
is $(\D C_m,\dots,\D C_1)$. In particular, we have $B = \D\D B = \D C_m.$

On the other hand, given a non-simple brick $M = B$, then $B$ is the only socle brick 
of $M$, but there are foundation bricks of $M$ which are proper submodules of $M$, thus
these foundation bricks are not socle bricks. 
$\s$
	\bigskip\medskip
{\bf 4. Some algebras with finite brick chain complexity.}
	\medskip
{\bf 4.1.} {\it \it For any $n\ge 3$, there is a directed gentle algebras with $n$ 
simple modules, with brick chain complexity equal to $n$.}
	\medskip
Proof. We consider the algebra $A$ with vertices $1,2,\dots,n$ and two
arrows $i \leftleftarrows i\!+\!1$ for $1\le i < n$ (always labeled $\alpha$ and $\beta$).
As relations, we take all the words $\alpha\beta$ and $\beta\alpha$, 
thus $A$ is gentle and directed. 
$$
{\beginpicture
    \setcoordinatesystem units <1cm,1cm>
\multiput{} at 0 0  0 4 /
\put{$1$} at 0 0 
\put{$2$} at 0 1
\put{$3$} at 0 2 
\put{$\vdots$\strut} at 0 2.6
\put{$n\!-\!1$} at 0 3
\put{$n$} at 0 4
\circulararc -30 degrees from -.1 .2 center at 1 0.5
\arr{-.12 0.3}{-.1 0.2}
\circulararc 30 degrees from  .1 .2 center at -1 .5
\arr{.12 0.3}{.1 0.2}

\circulararc -30 degrees from -.1 1.2 center at 1 1.5
\arr{-.12 1.3}{-.1 1.2}
\circulararc 30 degrees from  .1 1.2 center at -1 1.5
\arr{.12 1.3}{.1 1.2}

\circulararc -30 degrees from -.1 3.2 center at 1 3.5
\arr{-.12 3.3}{-.1 3.2}
\circulararc 30 degrees from  .1 3.2 center at -1 3.5
\arr{.12 3.3}{.1 3.2}

\multiput{$\alpha$\strut} at -.4 .5   -.4 1.5  -.4 3.5 /
\multiput{$\beta$\strut} at .4 .5     .4 1.5    .4 3.5 /
\endpicture} 
$$ 
Since $A$ is a string algebra, 
the indecomposable modules are strings and bands.
(Actually, the modules we are going to consider satisfy the
additional relations $\alpha^3$ and $\beta^2$.) 

A module annihilated by $\beta$ will be said to be an {\it $\alpha$-module.} An
indecomposable $\alpha$-module is uniquely determined by its support,
and we denote the indecomposable $\alpha$-module with support the interval $[a,b]$
by $I[a,b].$ 

For integers $1\le a < b \le n$, we are interested in the string module $N(a,b)$ with support
$[a,b]$, given by the word $(\alpha\beta^{-1}\alpha)^{b-a}$.
The following picture depicts on the left the module $N(1,5)$
(the action of $\alpha$ goes south-west, the action of $\beta$ south-east):
$$
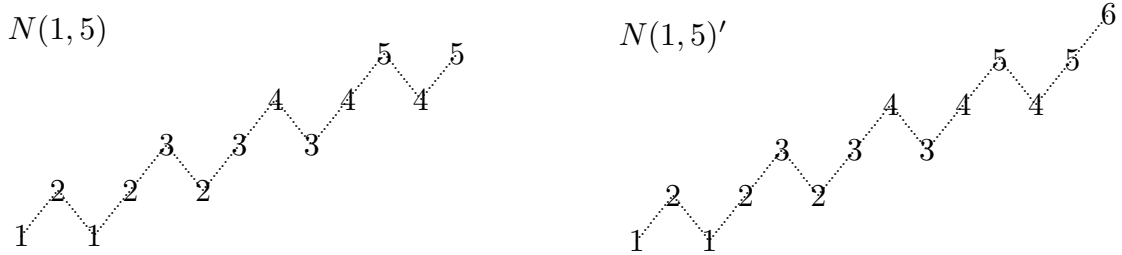
{\beginpicture
    \setcoordinatesystem units <.4cm,.5cm>
\put{\beginpicture
\multiput{} at 0 0  13 5 /
\multiput{$1$} at 0 0  2 0 /
\multiput{$2$} at 1 1  3 1  5 1 /
\multiput{$3$} at 4 2  6 2  8 2 /
\multiput{$4$} at 7 3  9 3  11 3 /
\multiput{$5$} at 10 4  12 4 /
\setdots <.5mm>
\plot 0 0  1 1  2 0  4 2  5 1  7 3  8 2  10 4  11 3  12 4 /
\put{$N(1,5)$} at 1 4.5 
\endpicture} at 0 0
\put{\beginpicture
\multiput{} at 0 0  13 5 /
\multiput{$1$} at 0 0  2 0 /
\multiput{$2$} at 1 1  3 1  5 1 /
\multiput{$3$} at 4 2  6 2  8 2 /
\multiput{$4$} at 7 3  9 3  11 3 /
\multiput{$5$} at 10 4  12 4 /
\multiput{$6$} at 13 5 /
\setdots <.5mm>
\plot 0 0  1 1  2 0  4 2  5 1  7 3  8 2  10 4  11 3  13 5 /
\put{$N(1,5)'$} at 1 4.5 
\endpicture} at 17 0
\endpicture}
$$
Our aim is to calculate the complexity of the modules $N(a,b)$. 
For the proof, we will need to look also at the string modules $N(a,b)'$ 
given by the words
$(\alpha\beta^{-1}\alpha)^{b-a}\alpha$, 
with integers $1\le a < b < n$, with support $[a,b\!+\!1].$ The
module $N(1,5)'$ is depicted on the right.

We will show: {\it The complexity of $N(a,b)$ is $b\!-\!a\!+\!1$, provided that $b\!-\!a\ge 2$,}
(since $N(a,a\!+\!1)$ is homogeneous, its complexity is $1$) and that 
{\it the complexity of $N(a,b)'$ is $b\!-\!a\!+\!1$, provided that $b\!-\!a\ge 1$.}

Looking at the pictures, let us stress the south-west lines: in both cases, there are
precisely 5 such lines, they correspond to indecomposable $\alpha$-modules and they
yield brick chain filtrations of length 5 (always, we start from the right:
the corresponding foundation brick for $N(1,5)$ is the $\alpha$-module $I[4,5]$;
for $N(1,5)'$, it is $I[4,6]$.
Of course, these brick chain filtrations
show that the complexity both of $N(a,b)$ and $N(a,b)'$ is at most $b\!-\!a\!+\!1.$
	
Using induction, we will prove for $N(a,b)$ and $N(a,b)'$ with $b\!-\!a\ge 2$, and also
for $N(a,a\!+\!1)'$, there are no brick chain filtrations of shorter length.
	\medskip
Proof. For the proof, we may assume that $a = 1$. 
First, let us consider $N(1,m)'$ with $1 < m < n.$ We show that the complexity of $N(1,m)'$
is at least $m$ (thus equal to $m$).

The induction starts with $m = 2:$  The module $M = N(1,2)'$ is obviously not
homogeneous, thus its complexity is at least $2$. 

Next, let $m \ge 3.$ Let $\Cal F$ be the set of the simple modules
$1,2,\dots,m\!-\!1$, and the $\alpha$-modules $I[m\!-\!1,m]$ and $I[m\!-\!1,m\!+\!1].$
We claim that {\it any foundation brick of $N(1,m)'$ belongs to $\Cal F.$}
For the proof, let $B$ be a foundation brick of $M$. 
We use 3.2 in order to specify possible foundation
bricks $B$. The simple modules occuring in the socle of $N(1,m)'$ are $1,2,\dots,m\!-\!1$,
thus either $B$ is one of the simple modules $1,2,\dots,m\!-\!1$, 
or else the top of $B$ is in $\add(m,m\!+\!1).$
The only submodule of $M$ with top $m\oplus(m\!+\!1)$ is not a brick, thus
the top of $B$ is either $m$ or $m\!+\!1$.
If the top of $B$ is $m\!+\!1$, then $B$ is the $\alpha$-module $I[m\!-\!1,m\!+\!1]$.
If the top of $B$ is $m$, and $B$ is either the $\alpha$-module $I[m\!-\!1,m]$ 
or else $B$ has length 4. The latter case is impossible, since in this case 
the trace of $B$ in $M$ has length 5, thus is not a direct sum of copies of $B$,
in contrast to 3.1. Altogether we see that $B$ belongs to $\Cal F.$

We are going to show now that for $B$ in $\Cal F$  
the factor module $M/t_B(M)$ has
complexity at least $m\!-\!1$. 

If $B = 1$, then $M/t_B(M)$ has $N(2,m)'$ as a direct summand.
If $B$ is equal to $m\!-\!1$ or $I[m\!-\!1,m]$ or $I[m\!-\!1,m\!+\!1],$ then 
$M/t_B(M)$ has $N(1,m\!-\!1)'$ as a direct summand. Again, in all these cases, 
$M/t_B(M)$ has complexity at least $m\!-\!1$. 

Now assume that $B = t$ with $2 \le t \le m\!-\!2.$ Then 
$M/t_B(M) = N(1,t)'\oplus N(t\!+\!1,m)'.$ A brick chain $\Cal B$ for $M/t_B(M)$ gives rise to
a brick chain $\Cal B'$ for $N(1,t)'$ and to a brick chain $\Cal B''$ for $N(t\!+\!1,m)'.$ 
By induction, the brick chain $\Cal B'$ has cardinality at least $t$, the brick chain
$\Cal B''$ has cardinality at least $m\!-\!t$. If $\Cal B'$ and $\Cal B''$ overlap, the
intersection can be only the simple module $t\!+\!1$, thus $\Cal B$ has cardinality
at least $t+(m\!-\!t)-1 = m\!-\!1,$ as we wanted to show. 
This completes the proof that $N(1,m)'$ has complexity $m$.
	\smallskip
Now let us consider $M = N(1,m)$ with $m\ge 3.$ The arguments are similar to those
used in the
previous induction step. This time, let $\Cal F$ be the set of the simple
modules $1,2,\dots,m\!-\!1$ and the $\alpha$-module $I[m\!-\!1,m].$ As above, we see that
any foundation brick of $N(1,m)$ belongs to $\Cal F$. We show that the
complexity of any factor module $M/t_B(M)$ with $B \in \Cal F$ is at least $m\!-\!1$. 

First, let $B = 1.$ If $m = 3$, one has to check that $M/t_B(M)$ is not homogeneous, thus
the complexity of $M/t_B(M)$ is at least $2$. If $m \ge 4,$ then $M/t_B(M)$ has
$N(2,m)$ as a direct summand, and, by induction, the complexity of $N(2,m)$ is at least
$m\!-\!1$. 

If $B$ is either $m\!-\!1$ or $I[m\!-\!1,m]$, then $M/t_B(M)$ has $N(1,m\!-\!1)'$ as a direct summand,
and we know already that the complexity of $N(1,m\!-\!1)'$ is $m\!-\!1.$ 

If $B = m\!-\!2,$ then $M/t_B(M) = N(1,m\!-\!2)\oplus N(m\!-\!1,m)$. 
Any brick chain $\Cal B$ for $M/t_B(M)$ gives rise to
a brick chain $\Cal B'$ for $N(1,m\!-\!2)'$.
As we know, the brick chain $\Cal B'$ has cardinality at least $m\!-\!2$. 
All the bricks in $\Cal B'$ have support in $[1,m\!-\!1]$, thus $\Cal B'$ is a proper
subsequence of $\Cal B$, therefore $\Cal B$ has cardinality at least $m\!-\!1.$

Finally, we have to consider the cases where $B = t$ with $2 \le t \le m\!-\!3.$ Then 
$M/t_B(M) = N(1,t)'\oplus N(t\!+\!1,m).$ 
A brick chain $\Cal B$ for $M/t_B(M)$ gives rise to
a brick chain $\Cal B'$ for $N(1,t)'$ and a brick chain $\Cal B''$ for $N(t\!+\!1,m).$ 
As we know, the brick chain $\Cal B'$ has cardinality at least $t$.
By induction,  the brick chain
$\Cal B''$ has cardinality at least $m\!-\!t$. If $\Cal B'$ and $\Cal B''$ overlap, the
intersection can be only the simple module $t\!+\!1$, thus $\Cal B$ has cardinality
at least $t+(m\!-\!t)-1 = m\!-\!1,$ as we wanted to show.
 
This completes the proof that $N(1,m)$ has complexity at least $m$.
$\s$
	\bigskip

{\bf 4.2.} The proof shows that {\it all bricks which occur in the brick type of a brick chain
filtration of $N(a,b)$ or of $N(a,b)'$ of minimal length 
are $\alpha$-modules}.

	\bigskip
{\bf 4.3.} In 4.1, we have assumed that $n \ge 3.$ What about $n = 2$ ?
{\it \it A directed gentle algebra with $2$ simple modules 
has brick chain complexity equal to $1$. There are algebras with $2$ simple modules which
are directed or gentle, which have brick chain complexity equal to $2$.}
	\medskip
Proof. The only directed gentle algebras with $2$ simple modules are the path algebras of the quivers with vertices $1,2$ and at most two arrows $1\leftarrow 2,$ thus factor algebras
of the Kronecker algebra. The Kronecker algebra has brick chain complexity 1.

The $t$-Kronecker algebra with $t\ge 3$ is directed and has brick chain complexity 2.
The Nakayama algebra with Kupisch series $(2,3)$ is gentle and has brick chain complexity 2.
$\s$
	\bigskip\bigskip
{\bf 5. The brick chain filtrations of a regular Kronecker module.}
	\medskip
As a preparation for section 6, we include some observations about Ktronecker modules;
these are the modules over the Kronecker algebra. The 
Kronecker algebra $A$ is the path algebra of the quiver with two vertices $1,2$
and two arrows $1\leftleftarrows 2$. 
There is the
well-known trisection of the indecomposable $A$-modules: there are the preprojective
modules, the regular modules $\Cal R$ and the preinjective modules.
The regular modules form an abelian category; non-isomorphic simple regular modules are
$\Hom$-orthogonal.
	\bigskip
{\bf 5.1. Lemma.} {\it 
Let $M$ be a regular Kronecker module. A foundation brick of $M$ is either
simple projective or else a simple regular module.}
	\medskip
Proof. Any indecomposable submodule $U$ of $M$ is preprojective or regular. 
If $U$ is preprojective and not simple, then $U$ generates $M$, thus the trace of $U$ in $M$
is not a direct sum of copies of $U$. According to 3.1, $U$ cannot be a foundation brick of
$M$. 
$\s$
	\bigskip
{\bf 5.2. Brick chain filtrations.}
Let $M$ be a non-zero regular Kronecker module with precisely $t$ isomorphism classes
of simple regular submodules. 
Then $M$ has the following brick chain filtrations:
	\medskip
There is one non-torsional brick chain filtration $(0 \subset M' \subset M)$ 
with $M' = \soc M = \rad M$, its brick type is $(1,2)$. 
If $B_1,\dots,B_t$ are the simple regular modules which
are submodules of $M$ (in any order), then $M = \bigoplus_i N_i$ with $N_i$ homogeneous of
type $B_i$. If $M_i = \bigoplus_{j\le i} N_j$, then $(M_i)_i$ is a torsional brick
chain filtration of brick type $(B_1,\dots,B_t).$ 
	\medskip
Thus. {\it $M$ has $t!$ torsional brick chain filtrations, all
of length $t$, and one non-torsional brick chain filtration of length $2$.}
	\bigskip\bigskip
{\bf 6. An algebra with infinite brick chain complexity.}
	\medskip
{\bf 6.1. The algebra.}
Let us present an algebra $A$ with brick chain complexity $\infty$. Take the path algebra
with two simple modules $1,2$, with two arrows $1\leftleftarrows 2$, two arrows 
$2\leftleftarrows 1$ and a loop at $1$; we assume that $1$ is a node (that means: all
paths $x \leftarrow 1 \leftarrow y$ are zero relations; there are five such paths), and that there
is no other relation. (Note that 
$A$ has Loewy length 3.  
After separation of the node, we deal with a hereditary algebra, namely the path algebra
of the quiver obtained from $1 \leftleftarrows 2 \leftleftarrows 1'$ by adding
an arrow $1\leftarrow 1'$).

Let $I'$ be the ideal generated by the two arrows $2 \leftleftarrows 1$ and the loop,
let $A' = A/I'.$ Similarly, let $I''$ be the ideal generated by the two arrows 
$1 \leftleftarrows 2$ and the loop, and let $A'' = A/I''.$ Then $A', A''$ are Kronecker algebras.
Here is the structure of the indecomposable projective $A$-modules $P_A(i),$ 
the indecomposable projective $A'$-modules $P_{A'}(i)$ and 
the indecomposable projective $A''$-modules $P_{A''}(i)$:
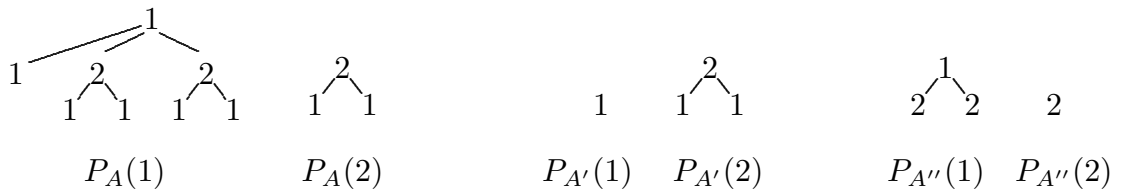
$$
{\beginpicture
    \setcoordinatesystem units <.3cm,.4cm>
\put{\beginpicture
\multiput{} at 0 0  2 2.5  /
\multiput{$1$} at 0 0  2 0 /
\multiput{$2$} at 1 1   /
\plot 0.3 0.3  0.7 0.7 /
\plot 1.3 0.7  1.7 0.3 /
\endpicture} at 8 0
\put{\beginpicture
\multiput{} at -1 0  2 2.5 /
\multiput{$1$} at -1 1  1 0  3 0  4 2.5  5 0  7 0 /
\multiput{$2$} at 2 1  6 1   /
\plot 3.6 2.3  -.5 1.3 /
\plot 3.7 2.1    2.3 1.6 /
\plot 4.3 2.1  5.7 1.6 /

\plot 1.3 0.3  1.7 0.7 /
\plot 2.3 0.7  2.7 0.3 /
\plot 5.3 0.3  5.7 0.7 /
\plot 6.3 0.7  6.7 0.3 /

\endpicture} at 0 0

\put{$P_A(1)$} at 0 -3
\put{$P_A(2)$} at 8 -3

\put{\beginpicture
\multiput{} at 0 0  2 2.5  /
\multiput{$1$} at 0 0  2 0  -3 0  /
\multiput{$2$} at 1 1  /
\plot 0.3 0.3  0.7 0.7 /
\plot 1.3 0.7  1.7 0.3 /
\endpicture} at 20 0

\put{\beginpicture
\multiput{} at 0 0  2 2.5  /
\multiput{$2$} at -4 0  -2 0  1 0  /
\multiput{$1$} at -3 1  /
\plot -2.7 0.7  -2.3 0.3 /
\plot -3.7 0.3  -3.3 0.7 /
\endpicture} at 32 0

\put{$P_{A'}(1)$} at 17 -3
\put{$P_{A'}(2)$} at 21.8 -3

\put{$P_{A''}(1)$} at 29.8 -3
\put{$P_{A''}(2)$} at 34.5 -3

\endpicture}
$$
The indecomposable $A'$-modules which are not simple have top in $\add 2$ and socle in 
$\add 1,$
The indecomposable $A''$-modules which are not simple have top in $\add 1$ and socle in 
$\add 2.$
	\bigskip
{\bf 6.2. The modules of interest.}
For natural numbers $t,u$ with 
$(t,u),$ we are going to exhibit $A$-modules $M = M(t,u)$ 
such that any brick chain filtration has length $t+1$.
Since we also show that for any $t\ge 1$, there are modules of the
form $M(t,0)$, which are indecomposable, 
we obtain in this way indecomposable $A$-modules with arbitrarily large
brick chain complexity.

Thus, let us describe the modules of the form $M(t,u).$ 
First, let us ignore the action of the loop. We take a direct sum $X$ of $t$ serial modules
$X_1,\dots, X_t$ of length three, each with top and socle equal to $1$ and with middle
composition factor equal to $2$, such that the modules $X'_i = \rad X_i$ are pairwise
non-isomorphic and such that also the modules $X''_i = X_i/\soc$ are pairwise
non-isomorphic. We denote by $N$ the socle of $X$.
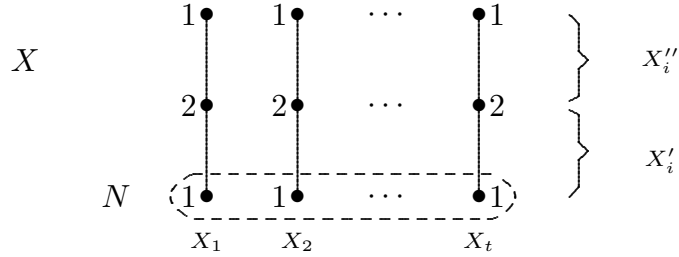
$$
{\beginpicture
    \setcoordinatesystem units <1cm,1cm>
\multiput{} at 0 0  3 2 /
\multiput{$\bullet$} at 0 0  0 2  1 0  1 2  3 0  3 2 
    0 1  1 1  3 1 /
\multiput{$1$} at -.2 0  -.2 2  0.8 0  .8 2  3.2 0  3.2 2 /
\multiput{$2$} at -.2 1  .8 1  3.2 1 /
\plot 0 0  0 2 /
\plot 1 0  1 2 /
\plot 3 0  3 2 /
\multiput{$\cdots$} at 2 0  2 1  2 2 /
\put{$\ss X_1$} at 0 -.5
\put{$\ss X_2$} at 1 -.5
\put{$\ss X_t$} at 3 -.5
\put{$\ss X''_i$} at  5 1.5
\put{$\ss X'_i$} at  5 0.4
\plot 4 0   4.1 0.1  4.1 0.4  4.2 0.5  4.1 0.6  4.1 0.9  4 0.95 /
\plot 4 1.05   4.1 1.1  4.1 1.4  4.2 1.5  4.1 1.6  4.1 1.9  4 2 /
\put{$X$} at -2 1.5
\put{$N$} at -1 0
\setquadratic
\setdashes <1mm>
\plot -.15 0.25  1.5 0.25  3.15 0.25  3.4 0  3.15 -.25  1.25 -.25  -.15 -.25  -.4 0  
   -.15 0.25 /
\endpicture}
$$
Next, let $Y_j$ be a copy of the simple module $1$, for $1\le j \le u$ and
$Y$ the direct sum of these modules.
Ignoring the action of the loop, the module $M(t,u)$ 
will be of the form $X\oplus Y$ and the action of the loop should send the top of $X$
into $N\oplus Y$ and should vanish on $(\rad X) \oplus Y$. Thus, 
the modules $M(t,u)$ can be visualized as follows
(the action of the loop is just indicated by some dashed line segments):

$$
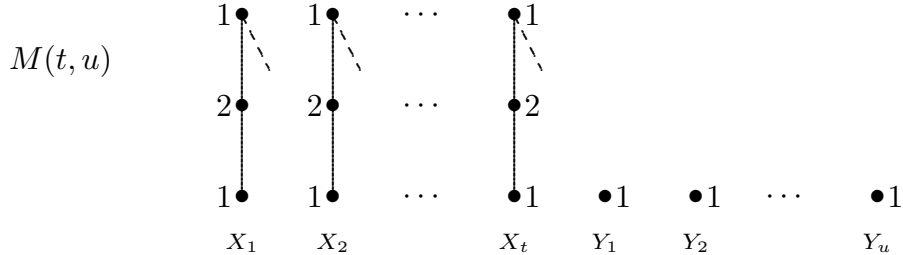
{\beginpicture
    \setcoordinatesystem units <1cm,1cm>
\multiput{} at 0 0  3 2 /
\multiput{$\bullet$} at 0 0  0 2  1 0  1 2  3 0  3 2 
    0 1  1 1  3 1 /
\multiput{$1$} at -.2 0  -.2 2  0.8 0  .8 2  3.2 0  3.2 2 /
\multiput{$2$} at -.2 1  .8 1  3.2 1 /
\plot 0 0  0 2 /
\plot 1 0  1 2 /
\plot 3 0  3 2 /

\setdashes <1mm>
\plot 0 2  0.3 1,4 /
\plot 1 2  1.3 1.4 /
\plot 3 2  3.3 1.4  /
\multiput{$\cdots$} at 2 0  2 1  2 2  6 0 /
\put{$\ss X_1$} at 0 -.5
\put{$\ss X_2$} at 1 -.5
\put{$\ss X_t$} at 3 -.5
\put{$\ss Y_1$} at 4 -.5
\put{$\ss Y_2$} at 5 -.5
\put{$\ss Y_u$} at 7 -.5

\multiput{$\bullet$} at 4 0  5 0  7 0 /
\multiput{$1$} at 4.2 0   5.2 0   7.2 0   /
\put{$M(t,u)$} at -2 1.5
\endpicture}
$$
	\bigskip
{\bf 6.3. Proposition.}  {\it If $M$ is a module of the form $M(t,u)$, then
any brick chain filtration of $M$ has length $t\!+\!1$.}
	\medskip
Proof, by induction on $t$. If $t = 0,$ then $M$ is a direct sum of copies of $1$, thus
there is just one brick chain filtraion, namely $(0 \subseteq M).$
We assume now that $t \ge 1.$

What are the foundation bricks of $M$? Of course, $1$ is a foundation brick.
The remaining ones are just the modules $\rad X_i.$ Namely, a foundation brick $B$ different from
$1$ must have top in $\add 2$, thus it is a submodule of $N = \rad M(t,u).$ 
This is an $A'$-module which is the direct sum of a regular Kronecker module and
copies of $1$, thus Lemma 5.1 shows that $B = X_i'$ for some $i$.

For $B = 1$, we have $M/t_B(M) = X'' = \bigoplus_i X''_i.$ This is a regular Kronecker module
for the algebra $A''$, thus we can use 5.2. Let $(M_i)_i$ be a brick chain filtraion of $X$
of brick type $(B_1,\dots,B_m)$ with $B_1 = 1.$ Then $M_1 = X'',$ and we deal with a
brick chain filtration of $M/M_1$ of brick type $(B_2,\dots,B_m).$ Since $B_1 = 1,$ 
and $\Hom(B_1,B_m) = 0,$ we see that $B_m \neq 1.$ 
According to 5.2, the brick type of any brick chain filtration of $X''$ is either $(2,1)$
or else consists of regular brick. As we have mentioned, the case $(2,1)$ is impossible,
thus we deal with a
brick chain filtration of $M/M_1$ of brick type $(B_2,\dots,B_m)$ with regular bricks
and the number $m-1$ of bricks is equal to $t.$ This shows that $m = t+1.$
	
Next, assume that $B = X_i'$ for some $1\le i \le t.$ Of course, we can assume that $i = t.$
The trace $t_B(M)$ of $B$ in $M$ is just $X'_t$, thus $M/t_B(M)$ is of the form $M(t-1,u+1)$,
thus, by induction any brick chain filtration of $M/t_B(M)$ has length $t$,
therefore the given filtration $(M_i)_i$ of $M$ with $M_1$ in $\Cal E(B)$ has length $t+1$.
Note that in the case $t = 1,$ we use that $M/t_B(M)$ is of the form $M(0,u+1)$, in
particular non-zero, so that $M_1$ is a proper submodule of $M$. $\s$
	\bigskip
{\bf 6.4.} {\it For any $m \ge 0,$ there are indecomposable $A$-modules with 
brick chain complexity equal to $m$.}
	\medskip
Proof. For $m = 1$, take a simple module. Now let $m\ge 2$. 
There are indecomposable modules of the
form $M(t,0)$ with $t = m-1.$ Namely, 
as action of the loop on $M(t,0)$, we 
send the top of $X_i$ onto the socle of $X_{i+1}$, for $1\le i < t.$
It is easy to see that then $M(t,0)$ is indecomposable.
According to 6.3, the module $M(t,0)$ has complexity $m$. 
$\s$
	\bigskip
{\bf 6.5. Remark.} A similar construction yields for $1\le s < t$ an indecomposable 
module $M$ which has
two kinds of brick chain filtrations, namely filtrations
of length $s+1$ and $t+1$ (and no other brick chain filtration):
As above, we first ignore the loop, and take the direct sum of $t$ serial modules
$X_1,\dots, X_t$ of length three, each with top and socle equal to $1$ and with middle
composition factor equal to $2$, such that the modules $X'_i = \rad X_i$ are pairwise
non-isomorphic, whereas now the modules $X''_i = X_i/\soc$ should belong to 
precisely $s$ isomorphism classes. Again, we use 
the action of the loop in order to achieve that $M$ is indecomposable. $\s$
	\bigskip\bigskip

{\bf 7. References.}
	\medskip
\item{[R1]} C\. M\. Ringel. Brick chain filtrations. A report. arXiv:2411.18427  
\item{[R2]} C\. M\. Ringel. In preparation.
    	\bigskip
    	\medskip
{\baselineskip=1pt
\rmk
Claus Michael Ringel\par
Fakult\"at f\"ur Mathematik, Universit\"at Bielefeld \par
POBox 100131, D-33501 Bielefeld, Germany  \par
ringel\@math.uni-bielefeld.de}
	\bigskip\bigskip

\bye